\documentclass[11pt,a4paper,reqno]{amsart}
\usepackage[T1]{fontenc}
\usepackage{lmodern}
\usepackage[a4paper,hmargin=24mm,top=26mm,bottom=27mm,includehead,
            headheight=14pt,headsep=12pt]{geometry}
\usepackage{amssymb,mathtools,microtype}
\usepackage[hidelinks]{hyperref}
\hypersetup{pdftitle={Isometric Classification of Sobolev Spaces}}
\numberwithin{equation}{section}
\newtheorem{theorem}{Theorem}[section]
\newtheorem{lemma}[theorem]{Lemma}
\newtheorem{fact}[theorem]{Fact}

\newcommand{\R}{\mathbb R}
\newcommand{\N}{\mathbb N}
\newcommand{\Q}{\mathbb Q}
\newcommand{\A}{\mathcal A}
\newcommand{\suppess}{\operatorname{ess\,supp}}
\newcommand{\id}{\operatorname{id}}
\newcommand{\norm}[1]{\left\lVert#1\right\rVert}
\newcommand{\one}{\mathbf 1}
\newcommand{\aee}{\text{a.e.}}
\allowdisplaybreaks[2]
\title[Isometries of real Sobolev spaces]{Isometric Classification of Sobolev Spaces}
\author[Samir Hamad]{Samir Hamad}
\date{5 November 2026}
\begin{document}
\begin{abstract}
We show that the norm obtained by summing the $L^1$ norms of all weak
derivatives through order $k$, including the function itself, determines
both $n$ and $k$ on the real Sobolev space $W^{k,1}(\R^n)$. Every surjective linear isometry
for this norm is given by a translation, a signed permutation of the
coordinates and multiplication by a fixed sign $c\in\{-1,1\}$.
\end{abstract}
\maketitle

\section{Introduction}

Throughout, functions are real-valued and operators are real-linear.
For integers $n,k\geq1$, we equip $W^{k,1}(\R^n)$ with the norm
\begin{equation}\label{eq:norm}
 N_{n,k}(f)=\sum_{\alpha\in\A_{n,k}}\norm{D^\alpha f}_1,
 \qquad
 \A_{n,k}=\{\alpha\in\N_0^n:|\alpha|\leq k\}.
\end{equation}
Here $|\alpha|=\alpha_1+\cdots+\alpha_n$, and $D^0f=f$. In particular, for $k=1$ the derivative
part is $\sum_j\|\partial_j f\|_1$. A signed permutation of the axes
permutes these terms. We shall see that the geometry of the norm
recovers the axes even when neither the dimension nor the order of the
target space is fixed in advance.

\begin{theorem}\label{thm:main}
Suppose that $n,m,k,s$ are positive integers and
\[
 T:(W^{k,1}(\R^n),N_{n,k})\longrightarrow
   (W^{s,1}(\R^m),N_{m,s})
\]
is a surjective linear isometry. Then $m=n$ and $s=k$. There
exist a signed permutation matrix $P$, a vector $b\in\R^n$ and
a sign $c\in\{-1,1\}$ such that, for every
$f\in W^{k,1}(\R^n)$,
\begin{equation}\label{eq:rigid}
 (Tf)(y)=c f(Py+b)\qquad\text{for almost every }y\in\R^n.
\end{equation}

\end{theorem}

A related question is whether, for $n,m\ge 2$ and $k,s\ge 1$,
a Banach-space isomorphism
\[
W^{k,1}(\mathbb{R}^n) \cong W^{s,1}(\mathbb{R}^m)
\]
necessarily implies $n=m$ and $k=s$.
For background on the isomorphic classification of Sobolev spaces,
see Pe{\l}czy{\'n}ski and Wojciechowski
\cite{PelczynskiWojciechowski}.~Diestel~\cite[Theorem~1]{Diestel} used an $L^p$ extension theorem
to classify surjective linear isometries of $W^{k,p}(U)$ on bounded
connected Lipschitz domains $U\subset\mathbb{R}^n$ in the sense of $(1.2)$. The result applies
for $1\leq p<\infty$ with $p\neq 2$, provided that the isometry maps
$W^{k,p}(U)\cap C(\overline U)$ onto itself. When $k>n/p$, every
function in $W^{k,p}(U)$ has a continuous representative on
$\overline U$ by Sobolev embedding, so this condition is automatic
\cite[Corollary~1]{Diestel}.
M.~Biegert and R.~Nittka~\cite{BiegertNittka} proved that, for
$1\leq p<\infty$ with $p\neq 2$, every linear isometry
$T:W^{1,p}(\Omega_1)\to W^{1,p}(\Omega_2)$ preserves disjointness.
Here $\Omega_1,\Omega_2\subset\mathbb{R}^n$ are open and the Sobolev
norms use the Euclidean length of the gradient. If both sets are
bounded and $T$ is order bounded as a map from $W^{1,p}(\Omega_1)$
into $L^p(\Omega_2)$, they also proved that there exist a measurable
map $\xi:\Omega_2\to\Omega_1$ and a real-valued weight $g$ on
$\Omega_2$ such that, for each $u\in W^{1,p}(\Omega_1)$,
\[
    (Tu)(y)=g(y)\,u(\xi(y))
    \qquad\text{for almost every }y\in\Omega_2.
\]
Taking $u=\mathbf{1}$ gives $g=T\mathbf{1}$ almost everywhere.

\section{Preliminaries}

We use the Lebesgue measure on Euclidean spaces and identify functions
modulo null sets. The support of a measurable function will always mean
the class of its nonzero set. Two functions are disjoint if their
product vanishes almost everywhere.

\begin{fact} {\cite[Theorem~6.19]{LiebLoss}}\label{lem:locality}
If $u\in W^{1,1}_{\mathrm{loc}}(\R^n;\R)$, then $\nabla u=0$ almost
everywhere on $\{u=0\}$. Consequently, for $f\in W^{k,1}(\R^n)$
every $D^\alpha f$ with $|\alpha|\leq k$ vanishes almost everywhere
on $\{f=0\}$.
\end{fact}

\begin{lemma}\label{lem:disjoint}
For $f,g\in W^{k,1}(\R^n)$ one has
\begin{equation}\label{eq:disjoint}
 fg=0\ \aee\quad\Longleftrightarrow\quad
 N_{n,k}(f+g)+N_{n,k}(f-g)=2N_{n,k}(f)+2N_{n,k}(g).
\end{equation}
In particular $T$ and $T^{-1}$ preserve disjointness.
\end{lemma}

To see this, suppose that $fg=0$ almost everywhere. By Fact~\ref{lem:locality},
$D^\alpha f$ vanishes almost everywhere on $\{f=0\}$, and $D^\alpha g$
vanishes almost everywhere on $\{g=0\}$. Since these two sets cover
$\R^n$ up to a null set, $D^\alpha f\,D^\alpha g=0$ almost everywhere
for every $|\alpha|\leq k$. Lamperti's equality criterion
\cite[Corollary~2.1]{Lamperti} therefore gives
\[
 \|D^\alpha f+D^\alpha g\|_1+\|D^\alpha f-D^\alpha g\|_1
 =2\|D^\alpha f\|_1+2\|D^\alpha g\|_1.
\]
Summing over $\alpha$ gives the norm identity in \eqref{eq:disjoint}.

Conversely, the triangle inequality bounds each term on the left by
the corresponding term on the right. Equality in the sum thus implies
equality for $\alpha=0$. Since $D^0f=f$ and $D^0g=g$, Lamperti's
criterion gives $fg=0$ almost everywhere. Both $T$ and $T^{-1}$
preserve the norm identity, so both preserve disjointness.

\section{The Sobolev jet}\label{sec:jet}

We use the product of Lebesgue and counting measure on
$\R^n\times\A_{n,k}$. A sheet is a set $\R^n\times\{\alpha\}$ with its
index fixed. We use $\gamma\leq\alpha$ for coordinatewise comparison,
$\alpha!=\prod_i\alpha_i!$ and
$\binom{\alpha}{\gamma}=\prod_i\binom{\alpha_i}{\gamma_i}$. Our approach of embedding Sobolev spaces into scalar $L^p$ spaces,
extending the induced isometry, and applying Lamperti's theorem
follows the strategy of Diestel \cite[Sections~2--3]{Diestel}.

\begin{lemma}\label{lem:jet-full}
Give $\A_{n,k}$ counting measure. The map
\[
 J_{n,k}f(x,\alpha)=D^\alpha f(x),\qquad
 J_{n,k}:W^{k,1}(\R^n)\longrightarrow
     L^1(\R^n\times\A_{n,k}),
\]
is an isometry onto a closed subspace. For $g(x)=e^{-|x|^2}$, the jet
$J_{n,k}g$ is nonzero almost everywhere. The ratios
$J_{n,k}f/J_{n,k}g$, as $f$ ranges over $W^{k,1}(\R^n)$, generate
the full product $\sigma$-algebra of $\mathbb{R}^n \times \mathcal{A}_{n,k}$ modulo null sets.
\end{lemma}

\begin{proof}
Integrating $|J_{n,k}f|$ over the sheet $\R^n\times\{\alpha\}$ gives
$\|D^\alpha f\|_1$. Summing over the sheets gives precisely $N_{n,k}(f)$,  i.e., $J_{n,k}$ is an isometry. 
Closedness now follows from completeness. We must find one jet which is
nonzero almost everywhere on every sheet, and show that ratios by this jet
distinguish the measurable sets of the whole product.

Let $g(x)=e^{-|x|^2}$. We show by induction on $|\alpha|$ that
\[
 D^\alpha g(x)=H_\alpha(x)\,g(x),
\]
where $H_\alpha$ is a polynomial whose unique term of highest total
degree is $(-2)^{|\alpha|}x^\alpha$. For $\alpha=0$ this holds with
$H_0=1$, since $D^0g=g$.

Suppose the assertion holds for $\alpha$ and write $d=|\alpha|$.
Since $\partial_j g=-2x_jg$, the product rule gives
\[
 D^{\alpha+e_j}g
 =\partial_j(H_\alpha g)
 =(\partial_jH_\alpha)g+H_\alpha(-2x_jg)
 =\bigl(\partial_jH_\alpha-2x_jH_\alpha\bigr)g.
\]
Thus $H_{\alpha+e_j}=\partial_jH_\alpha-2x_jH_\alpha$ is again a
polynomial. To identify its highest-degree term, note that
$\partial_jH_\alpha$ contains no term of degree $d+1$.
In $-2x_jH_\alpha$, the only term of degree $d+1$ comes from
multiplying the highest-degree term of $H_\alpha$ by $-2x_j$.
This gives
\[
 (-2x_j)(-2)^d x^\alpha=(-2)^{d+1}x^{\alpha+e_j}.
\]
All other terms have degree at most $d$, so this term cannot cancel.
This proves the induction step. In particular, $H_\alpha$ is a nonzero
polynomial of degree $|\alpha|$ for every $\alpha$. Remove its zero set for each of the
finitely many $\alpha$, and let $G$ be the common remainder. Since each $H_\alpha$ is a nonzero polynomial, $G$ has full
Lebesgue measure. On $G$ the jet of $g$ supplies a nonzero denominator on
every sheet.

For $a\in\Q^n$ put $g_a(x)=e^{a\cdot x}g(x)$. It is easy to see that $g_a\in
W^{k,1}$. On $G$ the product rule gives
\begin{equation}\label{eq:gaussian-profile}
 \frac{D^\alpha g_a(x)}{D^\alpha g(x)}
   =e^{a\cdot x}P_{x,\alpha}(a),\qquad
 P_{x,\alpha}(a)=
   \sum_{\gamma\leq\alpha}\binom{\alpha}{\gamma}
    a^{\alpha-\gamma}\frac{D^\gamma g(x)}{D^\alpha g(x)}.
\end{equation}
Suppose two pairs $(x,\alpha)$ and $(x',\alpha')$ give the same
quotients for all $a\in\Q^n$. The expressions are continuous in $a$, so
they agree for all $a\in\R^n$. If $x\ne x'$, set $v=x-x'$ and substitute $a=tv$ in the equality of the two
quotients. By \eqref{eq:gaussian-profile}, this gives
\[
 e^{t v\cdot x} P_{x,\alpha}(tv)
 = e^{t v\cdot x'} P_{x',\alpha'}(tv)
 \qquad (t\in\R).
\]
Dividing both sides by $e^{t v\cdot x'}$ and using
$v\cdot x-v\cdot x'=v\cdot(x-x')=v\cdot v=|v|^2$, we obtain
\[
 e^{t|v|^2} P_{x,\alpha}(tv)=P_{x',\alpha'}(tv)
 \qquad (t\in\R).
\]
Both sides have a polynomial factor taking value $1$ at $t=0$.
As $t\to+\infty$, the left side grows in absolute value faster than
any polynomial, a contradiction. Thus the quotients determine $x$.
With $x=x'$ the polynomial factors themselves agree. The unique term
of highest total degree in $P_{x,\alpha}(a)$ is
$a^\alpha g(x)/D^\alpha g(x)$, with nonzero coefficient. Comparing its
multiindex with that on the other sheet yields $\alpha=\alpha'$.

The quotients in \eqref{eq:gaussian-profile} therefore define an
injective Borel map
\[
\Phi:X:=G\times\A_{n,k}\longrightarrow\R^{\Q^n},
\qquad
\Phi(x,\alpha)
=
\left(
\frac{D^\alpha g_a(x)}{D^\alpha g(x)}
\right)_{a\in\Q^n}.
\]
By the Lusin--Souslin theorem in its standard Borel form
\cite[Corollary~15.2]{Kechris}, the image $\Phi(X)$ is Borel in
$\R^{\Q^n}$, and $\Phi$ is a Borel isomorphism from $X$ onto
$\Phi(X)$.

Let
\[
\mathcal S
:=
\sigma\!\left(
(x,\alpha)\longmapsto
\frac{D^\alpha g_a(x)}{D^\alpha g(x)}
:\ a\in\Q^n
\right).
\]
We claim that $\mathcal S=\mathcal B(X)$, where $\mathcal B(X)$
denotes the Borel $\sigma$-algebra of $X$.
Since each quotient is Borel measurable, we immediately have
$\mathcal S\subseteq\mathcal B(X)$.

For the reverse inclusion, let $E\in\mathcal B(X)$.
The Borel isomorphism property, together with the fact that
$\Phi(X)$ is Borel, implies that $\Phi(E)$ is Borel in
$\R^{\Q^n}$. Moreover, injectivity gives
\[
E=\Phi^{-1}\bigl(\Phi(E)\bigr).
\]
Since $\Q^n$ is countable, the Borel $\sigma$-algebra of
$\R^{\Q^n}$ is generated by the coordinate projections
$\pi_a$, $a\in\Q^n$. For each such projection,
\[
(\pi_a\circ\Phi)(x,\alpha)
=
\frac{D^\alpha g_a(x)}{D^\alpha g(x)}.
\]
Consequently, the inverse image under $\Phi$ of every Borel
subset of $\R^{\Q^n}$ belongs to $\mathcal S$.
In particular, $E=\Phi^{-1}(\Phi(E))\in\mathcal S$.
Thus $\mathcal B(X)\subseteq\mathcal S$, proving that the
coordinate quotients generate the Borel $\sigma$-algebra on
$G\times\A_{n,k}$. The omitted part of the product
is null. This proves the claimed generation modulo null sets.
\end{proof}

By Lemma~\ref{lem:jet-full} and Hardin's extension theorem
\cite[Theorem~4.2]{Hardin}, the isometry $J_{m,s}T J_{n,k}^{-1}$ between
the jet ranges extends to a surjective $L^1$ isometry $U$ of the ambient
spaces satisfying
\begin{equation}\label{eq:intertwine}
                    U J_{n,k}=J_{m,s}T .
\end{equation}

Put $X=\R^n\times\A_{n,k}$ and $Y=\R^m\times\A_{m,s}$, and denote
by $\mu$ and $\nu$ their respective products of Lebesgue measure
and counting measure. Write $\mathfrak M_X$ for the algebra of
measurable subsets of $X$ modulo $\mu$-null sets. Thus $[E]=[F]$
in $\mathfrak M_X$ precisely when $\mu(E\triangle F)=0$. Define
$\mathfrak M_Y$ in the same way using $\nu$.

By Lamperti's representation theorem
\cite[Theorem~3.1]{Lamperti}, the surjective linear isometry
$U:L^1(X,\mu)\to L^1(Y,\nu)$ admits a representation
\[
                     Uu=h\,\widehat{\tau}u,
\]
where
\[
             \tau:\mathfrak M_X\longrightarrow\mathfrak M_Y
\]
is a Boolean $\sigma$-isomorphism and $\widehat{\tau}$ is the
induced transformation on measurable functions, characterized by
\[
             \widehat{\tau}\mathbf 1_E=\mathbf 1_F
             \quad\text{whenever }[F]=\tau([E]).
\]

Surjectivity implies that $h\ne0$ $\nu$-almost everywhere.
Indeed, every function in the range of $U$ vanishes almost
everywhere on $Z:=\{h=0\}$. If $\nu(Z)>0$, $\sigma$-finiteness
provides a measurable set $B\subseteq Z$ with
$0<\nu(B)<\infty$. Then $\mathbf 1_B\in L^1(Y,\nu)$ cannot
belong to the range of $U$, contradicting surjectivity.

By construction, the induced transformation transports
measurable supports according to $\tau$:
\[
       \suppess(\widehat{\tau}u)=\tau(\suppess u),
       \qquad \suppess u:=[\{u\ne0\}].
\]
Since multiplication by $h$ does not change the measurable
support, we conclude that
\[
       \suppess(Uu)=\tau(\suppess u)
       \qquad (u\in L^1(X,\mu)).
\]
\section{Recovering the spatial algebra and the finite fibers}\label{sec:spatial}

Let $(E_i)_{i\in I}$ be a family of measurable subsets of
$(X,\mu)$. Its essential union is the least measurable upper
bound modulo null sets: it is represented by a measurable set
$E$ such that
\[
\mu(E_i\setminus E)=0 \qquad (i\in I),
\]
and whenever a measurable set $F$ satisfies
$\mu(E_i\setminus F)=0$ for every $i\in I$, one has
$\mu(E\setminus F)=0$.
Such an $E$ is unique modulo null sets.
On a $\sigma$-finite space, the essential union exists and is
represented by the union of a countable subfamily
\cite[Proposition~215B(iv)]{Fremlin}.

\begin{lemma}\label{lem:saturation}
Let $F\subset W^{s,1}(\R^m)$ be stable under multiplication by
$C_c^\infty(\R^m)$ functions. If $E$ is the essential union of the
supports of the functions in $F$, then $E\times\A_{m,s}$ is the
essential union of the supports of their jets.
\end{lemma}

\begin{proof}
For $f\in F$, Fact 2.1 gives $D^\beta f=0$ almost everywhere off $E$
for every $\beta$. This proves containment in $E\times\A_{m,s}$.
For the reverse containment, the difficulty is that $f(y)\ne0$ need not
imply $D^\beta f(y)\ne0$. We may, however, replace $f$ by an appropriate
smooth multiple without leaving $F$.

Choose $f_\ell\in F$ whose supports realize $E$, the case of empty $E$
being immediate. For $r\in\N$ take $\chi_r\in C_c^\infty(\R^m)$
equal to one on $B(0,r)$. For every $t\in\Q^m$ the multiplier
$y\mapsto e^{t\cdot y}\chi_r(y)$ is smooth with compact support. Fix $r,\ell\in\mathbb{N}$, $t\in\mathbb{Q}^m$, and
$\beta\in\mathcal{A}_{m,s}$. Since $\chi_r=1$ on $B(0,r)$,
the weak product rule gives
\[
D^\beta\!\left(e^{t\cdot(\cdot)}\chi_r f_\ell\right)(y)
=
e^{t\cdot y}
\sum_{\gamma\leq\beta}
\binom{\beta}{\gamma}\,
t^{\beta-\gamma}D^\gamma f_\ell(y),
\]
for almost every $y\in B(0,r)$. The exceptional set may depend on
$r,\ell,t,\beta$. Since there are only countably many such choices,
the union of all exceptional sets is a null set $N\subset\mathbb{R}^m$.
Thus, for each $y\notin N$, the identity holds for every
$\ell\in\mathbb{N}$, $t\in\mathbb{Q}^m$, and
$\beta\in\mathcal{A}_{m,s}$ whenever $r>|y|$.
If $f_\ell(y)\ne0$, the coefficient of $t^\beta$ in the sum is $f_\ell(y)$,
so the polynomial in $t$ is nonzero for some rational $t$. Choosing
$r>|y|$, we have found a member of $F$ whose jet is nonzero at
$(y,\beta)$. Since $r \in \mathbb{N}$ was arbitrary and almost every $y \in E \setminus N$ satisfies $f_p(y) \neq 0$ for some $p \in \mathbb{N}$, we obtain the conclusion.

\end{proof}

For a family $F\subset W^{k,1}(\R^n)$, define its disjoint complement
in $W^{k,1}(\R^n)$ by
\[
 F^d=
 \bigl\{
   u\in W^{k,1}(\R^n):
   uf=0\ \aee\ \text{for every }f\in F
 \bigr\}.
\]
For families in $W^{s,1}(\R^m)$, the disjoint complement is defined
analogously within $W^{s,1}(\R^m)$. For open $O\subset\R^n$ let $X_O$ consist of the Sobolev functions vanishing almost everywhere outside $O$. For an open ball $B$ one has
\begin{equation}\label{eq:ball-complement}
 X_B=(X_{\R^n\setminus\overline B})^d .
\end{equation}
Indeed, membership in $X_B$ implies disjointness from the exterior
family. Conversely, choose smooth bumps supported in the exterior whose
positive sets cover it, using a countable family of balls with closure
in the exterior. A function disjoint from these bumps vanishes almost
everywhere there. The only remaining set outside $B$ is $\partial B$,
which is null.

By surjectivity and Lemma~\ref{lem:disjoint}, $T(F^d)=(TF)^d$.
In particular $TX_B$ is a disjoint complement by (4.1). Such a complement is
stable under multiplication by smooth compactly supported functions,
since multiplying $u$ by a smooth factor preserves every relation $uf=0$.
On the source side, $X_B$ already has this
stability and its essential support union is $B$, as is seen using bumps
inside $B$. Write $E_B$ for the essential support union of $TX_B$.
Fix $f\in X_B$. Applying Lamperti's support formula
$\tau(\suppess u)=\suppess(Uu)$ to the function $u=J_{n,k}f$ gives
\[
  \tau\bigl(\suppess(J_{n,k}f)\bigr)
  =\suppess\bigl(U(J_{n,k}f)\bigr).
\]
By \eqref{eq:intertwine}, $U(J_{n,k}f)$ and $J_{m,s}(Tf)$ agree
almost everywhere on $\R^m\times\A_{m,s}$. Thus, for every $f\in X_B$,
\[
  \tau\bigl(\suppess(J_{n,k}f)\bigr)
  =\suppess\bigl(J_{m,s}(Tf)\bigr).
\]

We now take the essential union over $f\in X_B$. Since $\tau$ is an
order isomorphism of the measure algebras, it carries the least upper
bound of the supports on the left to the least upper bound of their
images. By Lemma~\ref{lem:saturation}, the essential union of the
supports of $J_{n,k}f$ for $f\in X_B$ is $B\times\A_{n,k}$.
As $f$ ranges over $X_B$, $Tf$ ranges over $TX_B$, so the same lemma
on the target gives $E_B\times\A_{m,s}$ as the essential union of
the supports of $J_{m,s}(Tf)$. Consequently,
\[
  \tau(B\times\A_{n,k})=E_B\times\A_{m,s}
  \qquad\text{modulo null sets}.
\]
We claim that the same formula holds for every Lebesgue measurable
set $C\subset\R^n$. Indeed, the sets $C$ for which there exists a
measurable $D\subset\R^m$ satisfying
\[
  \tau(C\times\A_{n,k})=D\times\A_{m,s}
  \qquad\text{modulo null sets}
\]
form a $\sigma$-algebra which contains the balls with centers in
$\mathbb{Q}^n$ and positive rational radii, so it contains every
Borel subset of $\R^n$. Every Lebesgue measurable set agrees with
a Borel set modulo a null set, and taking the product with the finite
set $\A_{n,k}$ preserves this agreement. This proves the claim.

Conversely, let $D\subset\R^m$ be Lebesgue measurable. Applying the
same argument to $T^{-1}$ and its extension on jets $U^{-1}$ gives
a measurable set $C\subset\R^n$ such that
\[
  \tau^{-1}(D\times\A_{m,s})=C\times\A_{n,k}
  \qquad\text{modulo null sets}.
\]
Thus every set of the form $D\times\A_{m,s}$ is the image under $\tau$
of a set of the form $C\times\A_{n,k}$, modulo null sets.

\begin{lemma}[The permutation of the derivative indices]\label{lem:fibers}
The isomorphism $\tau$ on the jet spaces acts by inverse images under a
map of the form
\[
 (y,\beta)\longmapsto(\varphi(y),\pi_y(\beta)).
\]
Here $\varphi:\R^m\to\R^n$ is a nonsingular isomorphism modulo null sets
and $\pi_y:\A_{m,s}\to\A_{n,k}$ is a bijection for almost every $y$.
\end{lemma}
\begin{proof}
The preceding argument gives an isomorphism between the Lebesgue
measure algebras of $\R^n$ and $\R^m$, determined by
\[
  \tau(C\times\A_{n,k})=D\times\A_{m,s}.
\]
By \cite[Theorem~343B(v) and Proposition~343F]{FremlinIII},
the isomorphism between the Lebesgue measure algebras of $\R^n$
and $\R^m$ induced by $\tau$ is represented by mutually inverse
nonsingular measurable maps
\[
  \varphi:\R^m\to\R^n,
  \qquad
  \psi:\R^n\to\R^m
\]
modulo null sets. In particular, for all Lebesgue measurable sets
$C\subset\R^n$ and $D\subset\R^m$,
\[
\begin{aligned}
  \tau(C\times\A_{n,k})
  &=\varphi^{-1}(C)\times\A_{m,s},\\
  \tau^{-1}(D\times\A_{m,s})
  &=\psi^{-1}(D)\times\A_{n,k},
\end{aligned}
\qquad\text{modulo null sets},
\]
and
\[
  \psi(\varphi(y))=y
  \quad\text{for almost every }y\in\R^m,
  \qquad
  \varphi(\psi(x))=x
  \quad\text{for almost every }x\in\R^n.
\]

For each $\alpha\in\A_{n,k}$, choose a measurable representative
$E_\alpha$ of $\tau(\R^n\times\{\alpha\})$. Since $\tau$ preserves
unions and intersections, the sets $E_\alpha$ are pairwise disjoint
modulo null sets and cover $\R^m\times\A_{m,s}$ modulo a null set.
There are only finitely many $\alpha$, so we may modify these sets
on a single null set to obtain an actual partition. Define
$\pi_y:\A_{m,s}\to\A_{n,k}$ by
\[
  \pi_y(\beta)=\alpha
  \quad\Longleftrightarrow\quad
  (y,\beta)\in E_\alpha.
\]
For every measurable $C\subset\R^n$, we have
\[
  C\times\{\alpha\}
  =(C\times\A_{n,k})\cap(\R^n\times\{\alpha\}),
\]
and hence
\[
  \tau(C\times\{\alpha\})
  =(\varphi^{-1}(C)\times\A_{m,s})\cap E_\alpha
  \qquad\text{modulo null sets}.
\]
The expression on the right is precisely the inverse image of
$C\times\{\alpha\}$ under the measurable map
\[
  F(y,\beta)=(\varphi(y),\pi_y(\beta)).
\]
Every measurable subset of $\R^n\times\A_{n,k}$ agrees modulo a
null set with a finite union of sets $C_\alpha\times\{\alpha\}$.
Thus $F$ represents $\tau$ by inverse images.

Applying the same construction to $\tau^{-1}$ gives a measurable
map
\[
  G(x,\alpha)=(\psi(x),\rho_x(\alpha)),
  \qquad
  \rho_x:\A_{n,k}\to\A_{m,s},
\]
which represents $\tau^{-1}$ by inverse images. Consequently,
$G\circ F$ and $F\circ G$ induce the identity on the respective
measure algebras. Each composition therefore equals the identity
almost everywhere. 
Indeed, since $G\circ F$ induces the identity on the measure algebra
of $\R^m\times\A_{m,s}$, for every measurable set
$D\subset\R^m\times\A_{m,s}$ the symmetric difference
\[
  N_D=(G\circ F)^{-1}(D)\mathbin{\triangle}D
\]
is a null set with respect to Lebesgue measure times counting measure.
Thus, for every $z\notin N_D$,
\[
  z\in D\quad\Longleftrightarrow\quad G(F(z))\in D.
\]
The exceptional set $N_D$ may depend on $D$. To obtain this equivalence
simultaneously for a family of sets that distinguishes every pair of
points, consider
\[
  \mathcal D=
  \bigl\{
    B(q,r)\times\{\beta\}:
    q\in\mathbb Q^m,\ r\in\mathbb Q,\ r>0,
    \ \beta\in\A_{m,s}
  \bigr\}.
\]
Here $B(q,r)$ denotes the open ball with center $q$ and radius $r$.
This family is countable, since the rational centers and radii form
countable sets and $\A_{m,s}$ is finite. It also separates points.

The countable union
\[
  N=\bigcup_{D\in\mathcal D}N_D
\]
is a null set. Consequently, for every $z\notin N$ and every
$D\in\mathcal D$,
\[
  z\in D\quad\Longleftrightarrow\quad G(F(z))\in D.
\]
If $G(F(z))\ne z$, the separating property of $\mathcal D$ gives
a set $D\in\mathcal D$ containing exactly one of these two points,
contradicting this equivalence. Hence $G(F(z))=z$ for every
$z\notin N$, which proves
\[
  G\circ F=\operatorname{id}_{\R^m\times\A_{m,s}}
  \qquad\text{almost everywhere}.
\]
Applying the same argument to $F\circ G$, using rational balls in
$\R^n$ and the finite index set $\A_{n,k}$, gives
\[
  F\circ G=\operatorname{id}_{\R^n\times\A_{n,k}}
  \qquad\text{almost everywhere}.
\]
Since $\A_{m,s}$ is finite, for almost every $y\in\R^m$ we have
\[
  G(F(y,\beta))=(y,\beta)
  \qquad\text{for every }\beta\in\A_{m,s}.
\]
Comparing the second coordinates gives
\[
  \rho_{\varphi(y)}\circ\pi_y=\id_{\A_{m,s}}.
\]
Similarly, for almost every $x\in\R^n$,
\[
  F(G(x,\alpha))=(x,\alpha)
  \qquad\text{for every }\alpha\in\A_{n,k},
\]
so
\[
  \pi_{\psi(x)}\circ\rho_x=\id_{\A_{n,k}}.
\]
By nonsingularity of $\varphi$, the latter identity remains valid
for almost every $y$ after substituting $x=\varphi(y)$.
Using $\psi(\varphi(y))=y$ almost everywhere, we obtain
\[
  \pi_y\circ\rho_{\varphi(y)}=\id_{\A_{n,k}}.
\]
Thus $\pi_y$ and $\rho_{\varphi(y)}$ are mutually inverse for almost
every $y$, which proves the asserted bijectivity.
\end{proof}

The weighted-composition representation of the ambient isometry can now
be written, for $u\in L^1(\R^n\times\A_{n,k})$, as
\[
 (Uu)(y,\beta)=h_\beta(y)u(\varphi(y),\pi_y(\beta)).
\]
For $a=J_{n,k}g$ from Lemma~\ref{lem:jet-full}, the weight on sheet
$\beta$ is
\[
 h_\beta(y)=\frac{(Ua)(y,\beta)}{a(\varphi(y),\pi_y(\beta))}.
\]
Both $a$ and $Ua$ have full support, so every $h_\beta$ is finite
and nonzero almost everywhere. Taking $u=J_{n,k}f$ and using \eqref{eq:intertwine} gives for almost all $y \in \mathbb{R}^m$:
\begin{equation}\label{eq:jet-formula}
 D^\beta(Tf)(y)=h_\beta(y)
         D^{\pi_y(\beta)}f(\varphi(y)),\qquad \beta\in\A_{m,s}.
\end{equation}
Also $|\A_{n,k}|=|\A_{m,s}|$. Appending the coordinate
$k-|\alpha|$ to $\alpha\in\A_{n,k}$ gives a bijection with the
$(n+1)$-tuples of nonnegative integers of sum $k$. The fiber count is
therefore
\begin{equation}\label{eq:fiber-count}
                   \binom{n+k}{k}=\binom{m+s}{s}.
\end{equation}
In \eqref{eq:jet-formula} the exceptional null set may depend on $f$.
We shall compare this identity only for countable families of tests.

\section{The zeroth and first derivative levels}\label{sec:levels}

So far the bijection $\pi_y$ may interchange derivatives of different
orders. Even the cardinality in \eqref{eq:fiber-count} does not determine
$n$ and $k$ separately. We first identify the index $0$, and then
identify the indices of length one.

\begin{lemma}\label{lem:zero-level}
One has $\pi_y(0)=0$ almost everywhere. In particular, if $h=h_0$, then
for every $f\in W^{k,1}(\R^n)$,
\begin{equation}\label{eq:composition}
 (Tf)(y)=h(y)f(\varphi(y))\qquad\text{for almost every }y.
\end{equation}
\end{lemma}
\begin{proof}
Take a ball $B$ with rational center and radius. Choose a cutoff
$\theta\in C_c^\infty(\R^n)$ equal to one on $B$. Consider
\[
 E=\varphi^{-1}(B)\cap\{y:\pi_y(0)\ne0\}.
\]
On $E$ the target value $T\theta$ reads a positive-order derivative
of $\theta$ by \eqref{eq:jet-formula}. Hence $T\theta=0$ almost
everywhere on $E$. By locality all its derivatives vanish there as well.

On the other hand, partition $E$ into the finitely many sets
$E_\beta=E\cap\{y:\pi_y(\beta)=0\}$. They cover $E$ up to a null set
because $\pi_y$ is onto. On $E_\beta$, the same jet formula gives
$D^\beta T\theta=h_\beta\,\theta\circ\varphi=h_\beta$.
Since $h_\beta$ is nonzero almost everywhere, $E_\beta$ is null. We conclude that $E$ is null. Choosing a cover of $\mathbb{R}^n$ by balls with rational centers and positive rational radii gives the conclusion.
\end{proof}

\begin{lemma}\label{lem:first-level}
For $1\leq j\leq m$, one has $|\pi_y(e_j)|=1$ almost everywhere.
Consequently $m=n$ and $s=k$.
\end{lemma}
\begin{proof}
Fix $j$, a rational ball $B\subset\R^n$, and a cutoff $\theta$ as in
the preceding proof. Applying \eqref{eq:composition} to $n+1$ functions
gives
\[
 v_0=T\theta=h,\qquad v_i=T(x_i\theta)=h\varphi_i
 \quad(1\leq i\leq n)
 \quad\text{a.e. on }\varphi^{-1}(B).
\]
All the $v_i$ belong to $W^{1,1}$. These identities allow us to compute
$\partial_j Tf$ locally without assuming weak differentiability of $h$ or
of $\varphi$.

Define the vector-valued function $v=(v_1,\ldots,v_n)$ and fix $\varepsilon>0$. We work on
$E_\varepsilon=\varphi^{-1}(B)\cap\{|h|\geq\varepsilon\}$.
Fix $f\in C_c^\infty(\R^n)$. Choose a smooth cutoff
$\eta_\varepsilon$ on $[0,\infty)$ equal to zero on $[0,\varepsilon/2]$
and to one on $[\varepsilon,\infty)$ with $0 \leq \eta_\varepsilon \leq 1$. For $(a,b)\in\R\times\R^n$ put
\[
 F_\varepsilon(a,b)=
 \begin{cases}
 \eta_\varepsilon(|a|)\,a f(b/a), & a\ne0,\\
 0, & a=0 .
 \end{cases}
\]
To apply the Sobolev chain
rule \cite[Theorem~4.2]{ArendtKreuter}, we show that $F_\varepsilon$ is globally Lipschitz by establishing
the estimate
\begin{equation}\tag{$*$}
  \left|DF_\varepsilon(a,b)[u,w]\right|
  \le C_\varepsilon\bigl(|u|+|w|\bigr)
\end{equation}
for all $(a,b)\in\mathbb R\times\mathbb R^n$ and all directions
$(u,w)\in\mathbb R\times\mathbb R^n$. The constant may depend on
$f$, $\varepsilon$ and the chosen cutoff, but is independent of
$a,b,u,w$. First fix $a\ne0$. To compute the derivative in the direction
$(u,w)$, consider the curve
\[
  t\longmapsto(a+tu,b+tw)
\]
and set
\[
  z(t)=\frac{b+tw}{a+tu},
  \qquad
  z=z(0)=\frac ba.
\]
The denominator is nonzero for all sufficiently small $t$.
Applying the quotient rule to each coordinate gives
\[
  z'(0)=\frac{aw-ub}{a^2}
       =\frac{w-uz}{a}.
\]
The product rule and the chain rule therefore yield
\[
  \begin{aligned}
    \left.\frac{d}{dt}\right|_{t=0}
    \left[(a+tu)f\!\left(\frac{b+tw}{a+tu}\right)\right]
    &=u f(z)+a\,\nabla f(z)\cdot z'(0)\\
    &=u f(z)+a\,\nabla f(z)\cdot\frac{w-uz}{a}\\
    &=u f(z)+\nabla f(z)\cdot w
      -u\,z\cdot\nabla f(z)\\
    &=\bigl(f(z)-z\cdot\nabla f(z)\bigr)u
      +\nabla f(z)\cdot w.
  \end{aligned}
\]

We also differentiate the factor $\eta_\varepsilon(|a|)$.
Since $a\ne0$,
\[
  \left.\frac{d}{dt}\right|_{t=0}|a+tu|
  =\operatorname{sgn}(a)\,u.
\]
Another application of the product rule gives
\[
  \begin{aligned}
    DF_\varepsilon(a,b)[u,w]
    ={}&\eta_\varepsilon(|a|)
    \Bigl[
      \bigl(f(z)-z\cdot\nabla f(z)\bigr)u
      +\nabla f(z)\cdot w
    \Bigr]\\
    &+\eta_\varepsilon'(|a|)
      \operatorname{sgn}(a)\,u\,a f(z).
  \end{aligned}
\]
Using $a\operatorname{sgn}(a)=|a|$, the last summand becomes
\[
  \eta_\varepsilon'(|a|)\,|a|\,f(z)\,u.
\]

To bound these expressions, choose $R>0$ such that
$\operatorname{supp}f\subset B(0,R)$ and put
\[
  M_0=\sup_{z\in\mathbb R^n}|f(z)|,
  \qquad
  M_1=\sup_{z\in\mathbb R^n}|\nabla f(z)|.
\]
Both constants are finite. Moreover, $\nabla f(z)=0$ outside
$B(0,R)$, so for every $z\in\mathbb R^n$,
\[
  |z\cdot\nabla f(z)|\le RM_1,
  \qquad
  |f(z)-z\cdot\nabla f(z)|\le M_0+RM_1.
\]
The derivative $\eta_\varepsilon'(|a|)$ can be nonzero only when
$\varepsilon/2\le|a|\le\varepsilon$. Consequently,
\[
  \left|
    \eta_\varepsilon'(|a|)\,|a|\,f(z)\,u
  \right|
  \le
  \varepsilon\|\eta_\varepsilon'\|_\infty M_0\,|u|.
\]
Combining these estimates and using $0\le\eta_\varepsilon\le1$,
we obtain
\[
  \begin{aligned}
    |DF_\varepsilon(a,b)[u,w]|
    \le{}&
    \bigl(
      M_0+RM_1+
      \varepsilon\|\eta_\varepsilon'\|_\infty M_0
    \bigr)|u|
    +M_1|w|.
  \end{aligned}
\]
This proves the desired estimate~$(*)$ for $a\ne0$.
For $|a|<\varepsilon/2$, the function $F_\varepsilon(a,b)$
vanishes for every $b\in\mathbb R^n$. Hence $F_\varepsilon$ is
also smooth in a neighborhood of every point $(0,b)$, with
derivative zero there. Thus~$(*)$ holds on all of
$\mathbb R^{n+1}$.

Integrating the bounded derivative along line segments shows that
$F_\varepsilon$ is globally Lipschitz. Since
$F_\varepsilon(0,0)=0$, the Sobolev chain rule
\cite[Theorem~4.2]{ArendtKreuter} gives
\[
  F_\varepsilon\circ(v_0,v)\in W^{1,1}(\mathbb R^m).
\]

For almost every $y\in E_\varepsilon$, we have
\[
\begin{aligned}
  \bigl(F_\varepsilon\circ(v_0,v)\bigr)(y)
  &=\eta_\varepsilon(|v_0(y)|)\,v_0(y)\,
    f\!\left(\frac{v(y)}{v_0(y)}\right)\\
  &=\eta_\varepsilon(|h(y)|)\,h(y)\,f(\varphi(y))\\
  &=h(y)\,f(\varphi(y))\\
  &=(Tf)(y).
\end{aligned}
\]
Here we used $v_0=h$, $v/v_0=\varphi$, and
$\eta_\varepsilon(|h|)=1$ on $E_\varepsilon$, since
$|h|\ge\varepsilon$ there.
Both $Tf$ and $F_\varepsilon\circ(v_0,v)$ belong to
$W^{1,1}(\mathbb R^m)$, and their difference vanishes almost
everywhere on $E_\varepsilon$. By locality, the first weak
derivatives of this difference also vanish almost everywhere
on $E_\varepsilon$. Hence
\[
  \partial_j(Tf)(y)
  =
  \partial_j\bigl(F_\varepsilon\circ(v_0,v)\bigr)(y)
  \qquad\text{for almost every }y\in E_\varepsilon.
\]
Combining this identity with the Sobolev chain rule gives
\[
  \partial_j(Tf)(y)
  =
  DF_\varepsilon\bigl(v_0(y),v(y)\bigr)
  \bigl[\partial_jv_0(y),\partial_jv(y)\bigr]
  \qquad\text{for almost every }y\in E_\varepsilon.
\]
On $E_\varepsilon$, we have
\[
  \eta_\varepsilon(|v_0(y)|)=1,
  \qquad
  \eta_\varepsilon'(|v_0(y)|)=0.
\]
Consequently, the formula for $DF_\varepsilon(a,b)[u,w]$ delivers
\[
\begin{aligned}
  \partial_j(Tf)(y)
  ={}&
  \bigl(
    f(\varphi(y))
    -\varphi(y)\cdot\nabla f(\varphi(y))
  \bigr)\partial_jv_0(y)\\
  &+\nabla f(\varphi(y))\cdot\partial_jv(y).
\end{aligned}
\]
Expanding the two scalar products yields
\[
\begin{aligned}
  \partial_j(Tf)(y)
  ={}&
  f(\varphi(y))\,\partial_jv_0(y)\\
  &-\sum_{i=1}^n
    \varphi_i(y)\,(\partial_i f)(\varphi(y))\,
    \partial_jv_0(y)\\
  &+\sum_{i=1}^n
    (\partial_i f)(\varphi(y))\,\partial_jv_i(y).
\end{aligned}
\]
Collecting the terms multiplying each
$(\partial_i f)(\varphi(y))$, we obtain
\[
\begin{aligned}
  \partial_j(Tf)(y)
  ={}&
  \partial_jv_0(y)\,f(\varphi(y))\\
  &+\sum_{i=1}^n
    \bigl(
      \partial_jv_i(y)-\varphi_i(y)\partial_jv_0(y)
    \bigr)
    (\partial_i f)(\varphi(y)).
\end{aligned}
\]

Define the coefficient functions on $E_\varepsilon$ by
\[
  a_0(y)=\partial_jv_0(y),
  \qquad
  a_i(y)=
  \partial_jv_i(y)-\varphi_i(y)\partial_jv_0(y),
  \quad 1\le i\le n.
\]
Then
\begin{equation}\label{eq:first-jet}
  \partial_j(Tf)(y)
  =
  a_0(y)f(\varphi(y))
  +\sum_{i=1}^n a_i(y)(\partial_i f)(\varphi(y))
  \qquad\text{a.e. on }E_\varepsilon.
\end{equation}

For the fixed ball $B$, cutoff $\theta$ and index $j$, the functions
$v_0=T\theta$ and $v_i=T(x_i\theta)$ are fixed independently of
$f$. The formulas defining $a_0,\ldots,a_n$ therefore give the same
coefficient functions for every $f\in C_c^\infty(\mathbb R^n)$.
For each such $f$, identity~\eqref{eq:first-jet} holds almost
everywhere on $E_\varepsilon$.

Let $\alpha\in\A_{n,k}$ have length at least two and consider
$E=E_\varepsilon\cap\{y:\pi_y(e_j)=\alpha\}$. On this set
\eqref{eq:jet-formula} identifies the left side of \eqref{eq:first-jet}
with $h_{e_j}(D^\alpha f)\circ\varphi$. Take successively the following
test functions.
\begin{enumerate}
\renewcommand{\labelenumi}{\textup{(\roman{enumi})}}
\item For $f=\theta$, $D^\alpha f$ and all first derivatives vanish on
      $B$, while $f=1$ there. Thus $a_0=0$ on $E$ almost everywhere.
\item For $f=x_i\theta$, the derivative $D^\alpha f$ still vanishes
      on $B$, whereas $\partial_\ell f=\delta_{i\ell}$ there. With
      $a_0=0$, the identity gives $a_i=0$ for every $i$.
\item For $f=x^\alpha\theta$, the right side of \eqref{eq:first-jet}
      is now zero on $E$, while $D^\alpha f=\alpha!$ on $B$.
      Hence $h_{e_j}\alpha!=0$ almost everywhere on $E$.
\end{enumerate}
Since $h_{e_j}$ is nonzero almost everywhere, $E$ has measure zero.
Only finitely many tests were needed for each fixed $B,j,\alpha$ and
$\varepsilon$. Take $\varepsilon=1/l$ for $l\in\N$ and let $B$ range
over the rational balls. The sets $E_\varepsilon$ obtained in this way cover
the target up to a null set because $h\ne0$ almost everywhere.
This excludes every $\alpha$ of length at least two. The value $0$ is
already used by $\pi_y(0)$, so injectivity gives $|\pi_y(e_j)|=1$.

The $m$ distinct first-order target indices consequently have distinct
images among the $n$ first-order source indices. Hence $m\leq n$.
Applying the same reasoning to $T^{-1}$ gives $n\leq m$. For fixed
$n\geq1$ the number $|\A_{n,k}|$ is strictly increasing in $k$,
since each higher order adds at least one multiindex. Thus
\eqref{eq:fiber-count} also gives $s=k$.
\end{proof}

Henceforth $m=n$ and $s=k$. There is a measurable permutation $\sigma_y$
of $\{1,\ldots,n\}$ at almost every $y$ such that
\begin{equation}\label{eq:first-permutation}
 \pi_y(e_j)=e_{\sigma_y(j)},\qquad
 \partial_j(Tf)(y)=h_j(y)(\partial_{\sigma_y(j)}f)(\varphi(y)),
 \qquad h_j=h_{e_j}.
\end{equation}

\section{Equal weights, truncation, and rigidity}\label{sec:rigidity}

\begin{lemma}\label{lem:weights}
For every $\beta\in\A_{n,k}$, one has $|h_\beta|=|h|$ almost
everywhere. For every nonnegative measurable function $v$,
\begin{equation}\label{eq:base-measure}
 \int_{\R^n}|h(y)|v(\varphi(y))\,dy=\int_{\R^n}v(x)\,dx .
\end{equation}
\end{lemma}
\begin{proof}
Fix $\alpha\in\mathcal A_{n,k}$ and a nonnegative function
$v\in L^1(\mathbb R^n)$. Define
$u\in L^1(\mathbb R^n\times\mathcal A_{n,k})$ by
\[
u(x,\gamma)=
\begin{cases}
v(x),& \gamma=\alpha,\\
0,& \gamma\ne\alpha.
\end{cases}
\]
Then
\[
\|u\|_{L^1(\mathbb R^n\times\mathcal A_{n,k})}
=\sum_{\gamma\in\mathcal A_{n,k}}
  \int_{\mathbb R^n}|u(x,\gamma)|\,dx
=\int_{\mathbb R^n}v(x)\,dx.
\]
For almost every $y\in\mathbb R^n$, the map
$\pi_y:\mathcal A_{n,k}\to\mathcal A_{n,k}$ is bijective.
There is therefore a unique index
$\beta_\alpha(y)\in\mathcal A_{n,k}$ such that
\[
\pi_y\bigl(\beta_\alpha(y)\bigr)=\alpha.
\]
The representation of $U$ consequently gives, for almost every
$y\in\mathbb R^m$ and every $\beta\in\mathcal A_{n,k}$,
\[
(Uu)(y,\beta)
=h_\beta(y)\,u\bigl(\varphi(y),\pi_y(\beta)\bigr)
=
\begin{cases}
h_\beta(y)\,v(\varphi(y)),& \beta=\beta_\alpha(y),\\
0,& \beta\ne\beta_\alpha(y).
\end{cases}
\]
The isometry of $U$ gives
\[
 \int_{\R^n}|h_{\beta_\alpha(y)}(y)|v(\varphi(y))\,dy=\int_{\R^n}v(x)\,dx.
\]
For $\alpha=0$ this is \eqref{eq:base-measure}, since
$\beta_0(y)=0$. Monotone convergence extends the formula from
integrable $v$ to all nonnegative $v$.

For any $R>0$, use as $v$ the indicator of a measurable subset of
$B(0,R)$. Comparing the formulas for $\alpha$ and $0$ shows that the
densities $|h_{\beta_\alpha(y)}(y)|$ and $|h(y)|$ have equal integrals
on each such preimage. Since $\varphi$ is invertible modulo null sets,
these preimages include, up to null sets, every measurable subset of
$\varphi^{-1}(B(0,R))$. Both densities have finite integral on this
set, so they agree almost everywhere there. Letting $R$ grow yields
the equality for a fixed $\alpha$, i.e., for each $\alpha\in\A_{n,k}$, we have a null set $N_\alpha\subset\R^m$ such that
\[
  |h_{\beta_\alpha(y)}(y)|=|h(y)|
  \qquad (y\notin N_\alpha).
\]
Since $\A_{n,k}$ is finite, the union
\[
  N:=\bigcup_{\alpha\in\A_{n,k}}N_\alpha
\]
is again a null set. Enlarge $N$ by the null set on which
$\pi_y$ fails to be a bijection. Thus, for every $y\notin N$,
the preceding equality holds simultaneously for all
$\alpha\in\A_{n,k}$.

Now fix $\beta\in\A_{m,s}$ and $y\notin N$, and set
$\alpha_y:=\pi_y(\beta)$. By the definition of the inverse indices,
\[
  \beta_{\alpha_y}(y)
  =\pi_y^{-1}\bigl(\pi_y(\beta)\bigr)
  =\beta.
\]
Since the equality is valid for every source index at this
point $y$, we may apply it with $\alpha=\alpha_y$ to obtain
\[
  |h_\beta(y)|
  =|h_{\beta_{\alpha_y}(y)}(y)|
  =|h(y)|.
\]
Consequently, $|h_\beta|=|h|$ almost everywhere for every
$\beta\in\A_{m,s}$, with a common exceptional null set.

\end{proof}

\begin{lemma}\label{lem:extend-W1}
The isometry $T$ extends to a surjective linear isometry of
$(W^{1,1}(\R^n),N_{n,1})$. The extension satisfies
\eqref{eq:composition} and \eqref{eq:first-permutation} with the same
base map and weights.
\end{lemma}
\begin{proof}
For $f\in W^{k,1}$, \eqref{eq:base-measure} gives
$\norm{Tf}_1=\norm f_1$. Moreover, summing \eqref{eq:first-permutation}
over $j$ and using $|h_j|=|h|$ and the fact that $\sigma_y$ is a permutation for almost all $y$, it follows that
\[
 \sum_{j=1}^n\norm{\partial_j Tf}_1
 =\int_{\R^n}|h(y)|\sum_{i=1}^n|(\partial_i f)(\varphi(y))|\,dy
 =\sum_{i=1}^n\norm{\partial_i f}_1,
\]
where the last equality follows from (6.1). Thus $T$ also preserves $N_{n,1}$ on $W^{k,1}$. This subspace is dense
in $W^{1,1}$ because it contains $C_c^\infty$. Extending $T$ and
$T^{-1}$ by continuity for the first-order norm gives inverse
isometries, as the compositions agree with the identity on $W^{k,1}$.
Write $\widetilde T$ for the continuous extension constructed above.
Fix $f\in W^{1,1}(\R^n)$ and choose
$f_\ell\in C_c^\infty(\R^n)$ such that
$f_\ell\to f$ in $W^{1,1}(\R^n)$.

By \eqref{eq:base-measure}, every $q\in L^1(\R^n)$ satisfies
\[
 \norm{h(q\circ\varphi)}_1=\norm{q}_1.
\]
Thus $q\mapsto h(q\circ\varphi)$ is a continuous linear map on
$L^1(\R^n)$. Moreover, for any tuple
$(q_i)_{i=1}^n$ with $q_i\in L^1(\R^n)$, the identities
$|h_j|=|h|$ and the fact that $\sigma_y$ is a permutation for
almost every $y$ give
\[
\begin{aligned}
 \sum_{j=1}^n\int_{\R^n}
 |h_j(y)|\,|q_{\sigma_y(j)}(\varphi(y))|\,dy
 &=\int_{\R^n}|h(y)|
   \sum_{i=1}^n|q_i(\varphi(y))|\,dy\\
 &=\sum_{i=1}^n\norm{q_i}_1,
\end{aligned}
\]
where the last equality follows by applying
\eqref{eq:base-measure} to each $|q_i|$.

Define
\[
 g(y)=h(y)f(\varphi(y)),
 \qquad
 G_j(y)=h_j(y)(\partial_{\sigma_y(j)}f)(\varphi(y)),
 \quad 1\le j\le n.
\]
The preceding identities show that $g,G_1,\ldots,G_n$ belong to
$L^1(\R^n)$. The value and derivative formulas already established
for $T$ apply to every $f_\ell$, since
$f_\ell\in W^{k,1}(\R^n)$. Applying the first norm identity to
$q=f_\ell-f$, we therefore obtain
\[
 \norm{Tf_\ell-g}_1
 =\norm{f_\ell-f}_1\longrightarrow0.
\]
Applying the identity for tuples to
$q_i=\partial_i f_\ell-\partial_i f$ gives
\[
 \sum_{j=1}^n\norm{\partial_j(Tf_\ell)-G_j}_1
 =\sum_{i=1}^n\norm{\partial_i f_\ell-\partial_i f}_1
 \longrightarrow0.
\]
In particular, $\partial_j(Tf_\ell)\to G_j$ in $L^1(\R^n)$
for each $j$.

On the other hand, the continuity of $\widetilde T$ on
$W^{1,1}(\R^n)$ and its agreement with $T$ on $W^{k,1}(\R^n)$
imply
\[
 Tf_\ell=\widetilde T f_\ell
 \longrightarrow\widetilde T f
 \qquad\text{in }W^{1,1}(\R^n).
\]
Consequently,
\[
 Tf_\ell\longrightarrow\widetilde T f,
 \qquad
 \partial_j(Tf_\ell)\longrightarrow\partial_j(\widetilde T f)
 \qquad\text{in }L^1(\R^n).
\]
Uniqueness of the $L^1$ limit now gives
$\widetilde T f=g$ and $\partial_j(\widetilde T f)=G_j$.
Hence
\[
 (\widetilde T f)(y)=h(y)f(\varphi(y))
\]
and
\[
 \partial_j(\widetilde T f)(y)
 =h_j(y)(\partial_{\sigma_y(j)}f)(\varphi(y)),
 \qquad 1\le j\le n,
\]
for almost every $y\in\R^n$. Thus both formulas remain valid
for the extension, with the same map $\varphi$, the same
permutations $\sigma_y$, and the same weights.

\end{proof}

We now use the first-order extension to truncate a preimage under $T$.
A Lipschitz truncation preserves $W^{1,1}$, whereas for $k>1$ there is
no such general rule on $W^{k,1}$.

\begin{lemma}[The weight is constant]\label{lem:constant}
There is a constant $c\in\R\setminus\{0\}$ such that $h=c$ almost everywhere.
\end{lemma}
\begin{proof}
We work with the first-order extension. Fix a rational ball $K \subset \mathbb{R}^n$ and a positive integer $a$. Choose $u\in C_c^\infty(\R^n)$
equal to $a$ on $K$. By surjectivity there exists $f\in W^{1,1}$ with
$Tf=u$. The value formula gives $f(\varphi(y))=a/h(y)$ on $K$.
Also $\partial_j u=0$ there for each $j$. Since all $h_j$ are nonzero
and $\sigma_y$ is a permutation, \eqref{eq:first-permutation} yields
$(\partial_i f)(\varphi(y))=0$ for all $i$ at almost every $y\in K$.

We now truncate $|f|$. The function
$\psi(t)=\min(|t|,1)$ is $1$-Lipschitz on $\R$ and $\psi(0)=0$.
The Sobolev composition estimate \cite[Theorem~4.2]{ArendtKreuter} gives $\psi(f)\in W^{1,1}$ and
$|\partial_i(\psi(f))|\leq|\partial_i f|$ almost everywhere.
The preceding inequalities hold for all $i=1,\ldots,n$ outside
a common Lebesgue null set $N\subset\R^n$. Since $\varphi$ is
nonsingular, the inverse image $\varphi^{-1}(N)$ is also a
Lebesgue null set. Thus, for almost every $y\in\R^n$, we have
$\varphi(y)\notin N$, and consequently
\[
 \left|(\partial_i(\psi(f)))(\varphi(y))\right|
 \le \left|(\partial_i f)(\varphi(y))\right|,
 \qquad 1\le i\le n.
\]
Applying \eqref{eq:first-permutation} to $\psi(f)$
therefore gives $\nabla T(\psi(f))=0$ almost everywhere on $K$, since $(\partial_i f)(\varphi(y))=0$ for all $i$ at almost every $y \in K$.
Its value on $K$ is, by \eqref{eq:composition},
\[
 T(\psi(f))(y)=h(y)\min\!\left(\frac{a}{|h(y)|},1\right)=:h_a(y).
\]

The expression $h_a$ is independent of the choice of $u$ and of $K$.
Since $\nabla T(\psi(f))=0$ almost everywhere on $K$, we have that $T(\psi(f))$ is almost everywhere constant there. For a fixed $a$, apply this
conclusion on the countable nested balls $B(0,R)$, $R\in\N$.
Their constants agree on overlaps, so $h_a=c_a$ almost everywhere
on $\R^n$ for one scalar $c_a$.

Finally take a common full-measure set for all integers $a\geq1$ on which
$0<|h|<\infty$. At every point of this set $h_a(y)=h(y)$ once
$a\geq|h(y)|$. Thus the same scalar sequence $(c_a)$ converges to
$h(y)$ at every such $y$. Its limit gives the required nonzero $c$.
\end{proof}

\begin{lemma}\label{lem:rigidity}
The base map has a representative $\varphi(y)=Py+b$, where $P$ is a
signed permutation matrix. The constant in Lemma~\ref{lem:constant}
belongs to $\{-1,1\}$.
\end{lemma}
\begin{proof}
Choose $\theta\in C_c^\infty(\R^n)$ equal to one on $B(0,1)$.
For integers $R\geq1$ and $1\leq i\leq n$, put
$f_{R,i}(x)=x_i\theta(x/R)$. Their derivative is
\[
 \partial_\ell f_{R,i}(x)
  =\delta_{i\ell}\theta(x/R)+(x_i/R)(\partial_\ell\theta)(x/R).
\]
With $z=x/R$, the second term is $z_i(\partial_\ell\theta)(z)$, which
has a bound independent of $R$. Put
$F_{R,i}=c^{-1}Tf_{R,i}=f_{R,i}\circ\varphi$.
Each $F_{R,i}$ is bounded and its weak first derivatives are bounded
uniformly in $R$ by \eqref{eq:first-permutation}, since $|h_j/c|=1$.
The $W^{1,\infty}$ representative theorem \cite[Theorem~4.1]{Heinonen} gives globally Lipschitz
representatives with a common Lipschitz bound.

Take a full-measure set on which all these representatives agree with
$f_{R,i}\circ\varphi$. For any $y$ in this set,
$F_{R,i}(y)=\varphi_i(y)$ whenever $R>|\varphi(y)|$. Passing to the
limit $R \rightarrow \infty$ in the Lipschitz inequality at two points of this set gives the
Lipschitz bound for $\varphi_i$. Extension from this dense set gives a
globally Lipschitz representative of each coordinate of $\varphi$.
On $\{y:|\varphi(y)|<R\}$ the functions $F_{R,i}$ and $\varphi_i$
agree. By locality their derivatives agree there almost everywhere, and
$\partial_\ell f_{R,i}=\delta_{i\ell}$ on $B(0,R)$.
Using \eqref{eq:first-permutation} and exhausting in $R$, we obtain
\begin{equation}\label{eq:Dphi}
 \partial_j\varphi_i(y)=\frac{h_j(y)}c\,
       \one_{\{i=\sigma_y(j)\}}\qquad\text{almost everywhere}.
\end{equation}

Each column of $D\varphi$ has one nonzero entry, and these entries occupy distinct rows, since $\sigma_y$ is a permutation for almost every $y$. Each nonzero entry is a ratio $h_j(y)/c$ of
real numbers with $|h_j(y)|=|c|$, so it is either $1$ or $-1$. Thus $D\varphi(y)$ is a
signed permutation matrix almost everywhere. Its Euclidean operator
norm is one and therefore
$|\varphi(y)-\varphi(z)|\leq|y-z|$.

Apply the same argument to $T^{-1}$, using $U^{-1}$ as its ambient
extension. This gives a globally $1$-Lipschitz representative $\psi$
of the inverse base map. Since both original maps are nonsingular,
replacing them by their continuous representatives preserves both
inverse identities almost everywhere. By continuity the identities
hold everywhere. Applying the bound for $\psi$ at $\varphi(y)$ and
$\varphi(z)$ gives the opposite inequality, so that
\[
 |y-z|=|\psi(\varphi(y))-\psi(\varphi(z))|
       \leq|\varphi(y)-\varphi(z)|\leq|y-z|.
\]
Thus $\varphi$ is a surjective Euclidean isometry. By the Mazur--Ulam
theorem \cite{MazurUlam}, $\varphi(y)=Py+b$. Its linear part $P$
preserves the Euclidean norm and is therefore orthogonal.
Its derivative is the constant matrix $P$, and \eqref{eq:Dphi}
therefore makes $P$ a signed permutation matrix. Since $h=c$ and $\varphi(y)=Py+b$, identity
\eqref{eq:base-measure} reads
\[
 |c|\int_{\R^n}v(Py+b)\,dy
 =\int_{\R^n}v(x)\,dx.
\]
Let $B=B(0,1)$ and take $v=\mathbf{1}_B$. Then
\[
 |c|\int_{\R^n}\mathbf{1}_B(Py+b)\,dy=|B|.
\]
The change of variables $x=Py+b$ gives
$dy=|\det P|^{-1}\,dx$, and hence
\[
\begin{aligned}
 |c|\int_{\R^n}\mathbf{1}_B(Py+b)\,dy
 &=|c|\,|\det P|^{-1}
   \int_{\R^n}\mathbf{1}_B(x)\,dx\\
 &=|c|\,|\det P|^{-1}|B|.
\end{aligned}
\]
Consequently,
\[
 |c|\,|\det P|^{-1}|B|=|B|.
\]
Since $0<|B|<\infty$, division by $|B|$ yields
\[
 |c|\,|\det P|^{-1}=1.
\]
The matrix $P$ is orthogonal, so $|\det P|=1$.
Thus $|c|=1$, and since $c$ is real, we conclude that
$c\in\{-1,1\}$.
\end{proof}

Formula~\eqref{eq:composition} now gives \eqref{eq:rigid}, completing
the proof of Theorem~\ref{thm:main}.

\section*{Acknowledgments}

The author used OpenAI's GPT-6 Astra to assist with mathematical
exploration, the development and checking of proofs, and the
presentation of arguments. The author takes full responsibility
for the mathematical content and any remaining errors.

\end{document}